\documentclass[a4paper,10pt]{article}

\usepackage{amssymb}
\usepackage{amsmath}
\usepackage{amsthm}

\makeatletter
\let\@afterindentfalse\@afterindenttrue
\@afterindenttrue
\makeatother

\newtheorem{definition}{Definition}
\newtheorem{theorem}{Theorem}
\newtheorem{corollary}{Corollary}
\theoremstyle{remark} 
\newtheorem{remark}{Remark}[section]

\DeclareMathOperator*{\smallwedge}{\wedge}
\DeclareMathOperator*{\bigosquare}{
  \settowidth{\unitlength}{$\square$}
  \begin{picture}(1,1)(0,0)
  \put(0.1,0){\line(1,0){0.8}}
  \put(0.1,0.8){\line(1,0){0.8}}
  \put(0.1,0){\line(0,1){0.8}}
  \put(0.9,0){\line(0,1){0.8}}
  \put(0.5,0.4){\circle*{0.3}}
  \end{picture}}
\newcommand{\sign}{\mathop{\mathrm{sgn}}\nolimits}
\newcommand{\id}{\mathop{\mathrm{id}}\nolimits}
\newcommand{\Lin}{\mathop{\mathrm{Lin}}\nolimits}
\newcommand{\btop}[2]{\genfrac{}{}{0pt}{2}{#1}{#2}}
\newcommand{\dint}{\mathop{\mathrm{d}}\nolimits}
\newcommand{\rank}{\mathop{\mathrm{rank}}\nolimits}

\begin{document}

\title{Short note on the perturbation of operators with dyadic products
\thanks{keywords: perturbation, dyadic product, inverse operator;  MSC: 15A09, 15A15}
  }
\author{Attila Andai\thanks{andaia@math.bme.hu}\\
  RIKEN, BSI, Amari Research Unit \\
  2--1, Hirosawa, Wako, Saitama 351-0198, Japan.}
\date{October 7, 2008}

\maketitle

\begin{abstract}
In this paper we use abstract vector spaces and their duals without any canonical basis.
Some of our results can be extended to infinite dimensional vector spaces too, but
  here we consider only finite dimensional spaces.
We focus on a general perturbation problem.
Assume that $B:V\to V$ is a linear operator, which is perturbated to $B'=B+Q$.
We examine the question how the determinant and the inverse change, because of
  this perturbation.
In our approach the operator $Q$ is given as a sum of dyadic products  
  $Q=\sum_{i=1}^{k}v_{i}\otimes p_{i}$, where  $v_{i}\in V$ and $p_{i}\in V^{*}$.
In this paper we derive an $m$-th order ($m\in\mathbb{N}$) approximation formula for 
  $\det B'$ and $(B')^{-1}$, which gives the exact result if $m\geq k$.
\end{abstract}

\section*{Introduction}

In this paper we use abstract vector spaces and their duals without any canonical basis.
Our notations follow the classical linear algebraic notations, for details see for example 
  \cite{Bha}.
We assume that $V$ is a real or complex vector space and the maps between vector spaces are linear.
Some of our results can be extended to infinite dimensional vector spaces too, but
  here we consider only finite dimensional spaces.

We focus on a general perturbation problem.
Assume that $B:V\to V$ is a linear operator, which is perturbated to $B'=B+Q$.
We examine the question how the determinant and the inverse change, because of
  this perturbation.
One natural approximation is given by the Taylor expansion, but that process requires
  norm on the vector space, and the Taylor polynomials do not give the exact result when 
  the series is cut within a finite number of terms.
In our approach the operator $Q$ is given as a sum of dyadic products 
\begin{equation*}
Q=\sum_{i=1}^{k}v_{i}\otimes p_{i},
\end{equation*}
  where  $v_{i}\in V$ and $p_{i}\in V^{*}$ for every $i=1,\dots,k$; 
  moreover we can assume that $k\leq\dim V$.
In this paper we derive an $m$-th order ($m\in\mathbb{N}$) approximation formula for 
  $\det B'$ and $(B')^{-1}$, which gives the exact result if $m\geq k$.

\section{On the inverse of perturbated operators}

If $A:V\times V^{*}\to\mathbb{R}$ is a bilinear map, then there exists a unique
  $\kappa(A):V\to V$ map such that
\begin{equation*}
 p(\kappa(A)x)=A(x,p) \quad \forall (x,p)\in V\times V^{*}.
\end{equation*}
We define the map $\kappa:\Lin(V\times V^{*},\mathbb{R})\to\Lin(V,V)$, which is
  an isomorphism.

\begin{definition}
Let $V$ be a vector space, such that $\dim V=n\geq 2$.
For given vectors $z_{1},\dots,z_{k}\in V$ and covectors $p_{1},\dots,p_{k}\in V^{*}$, where
  $1\leq k\leq n-1$ define the map
\begin{equation}
\Phi_{z_{1},\dots,z_{k}}^{p_{1},\dots,p_{k}}:V\times V^{*}\to\mathbb{R}\qquad
  (v,q)\mapsto \left(q\wedge p_{1}\wedge\dots\wedge p_{k} \right)(v,z_{1},\dots,z_{k}).
\end{equation}
We have $\kappa\left(\Phi_{z_{1},\dots,z_{k}}^{p_{1},\dots,p_{k}}\right)\in\Lin(V,V)$, and
  we introduce the symbol
\begin{align}
\bigosquare_{i=1}^{k}(z_{i},p_{i})=
  (z_{1},p_{1})\bigosquare(z_{2},p_{2})\bigosquare\dots\bigosquare(z_{k},p_{k})
  =\kappa\left(\Tilde{\Phi}_{z_{1},\dots,z_{k}}^{p_{1},\dots,p_{k}} \right).
\end{align}
\end{definition}
It is obvious from the definition that if $(z_{i})_{i=1,\dots,k}$ or $(p_{i})_{i=1,\dots,k}$ 
  are not linearly independent, then
\begin{equation*}
\bigosquare_{i=1}^{k}(z_{i},p_{i})=0
\end{equation*}
  and if $\pi,\pi'$ are permutations of the set $\left\{1,\dots,k\right\}$ then we have
\begin{equation*}
\bigosquare_{i=1}^{k}(z_{\pi(i)},p_{\pi'(i)})
  =(-1)^{\sign(\pi)\sign(\pi')}\bigosquare_{i=1}^{k}(z_{i},p_{i}).
\end{equation*}
Moreover if $v\in V$ and $q\in V^{*}$ then the equalities
\begin{align*}
&q\left[\left(\bigosquare_{i=1}^{k}(z_{i},p_{i})\right)(v)\right]=
  \Bigl(q\wedge p_{1}\wedge\dots\wedge p_{k}\Bigr) (v,z_{1},\dots,z_{k})\\
&\label{eq:osquare3}\left(\bigosquare_{i=1}^{k}(z_{i},p_{i})\right)(v)=
  v(p_{1}\wedge\dots\wedge p_{k})(z_{1},\dots,z_{k})
  -\sum_{i=1}^{k}z_{i}(p_{1}\wedge\dots\wedge p_{k})
  (z_{1},\dots,z_{i-1},v,z_{i+1},\dots,z_{k})
\end{align*}
hold.

Now we show how one can compute the determinant and the inverse of the perturbated
  identity operator.

\begin{theorem}
Let $V$ be a finite dimensional vector space, $\dim V=n$ and $(u_{i})_{i=1,\dots,k}$ a 
  family of vectors and $(p_{i})_{i=1,\dots,k}$ a family of covectors.
Define the linear map
\begin{equation}
A=\id_{V}+\sum_{i=1}^{k}u_{i}\otimes p_{i}.
\end{equation}
We have for the determinant of $A$ 
\begin{equation}
\label{eq:Adetmain} \det A=1+\sum_{i=1}^{\min(n,k)} \sum_{1\leq j_{1}<j_{2}<\dots<j_{i}\leq k} 
  \left(\smallwedge_{l=1}^{i}p_{j_{l}}\right)\left(u_{j_{1}},\dots,u_{j_{i}} \right)
\end{equation}
  and if $\det A\neq 0$, then we have for the inverse of $A$
\begin{equation}
\label{eq:Ainv}
A^{-1}\det A=\id_{V}+\sum_{i=1}^{\min(n-1,k)} 
  \sum_{1\leq j_{1}<j_{2}<\dots<j_{i}\leq k} \bigosquare_{l=1}^{i}(u_{j_{l}},p_{j_{l}}),
\end{equation}
  which can be written in the following form: if $x\in V$ and $q\in V^{*}$ then
\begin{equation}
\label{eq:Ainvmain}
q(A^{-1}x)\det A=q(x)+\sum_{i=1}^{\min(n-1,k)} \sum_{1\leq j_{1}<j_{2}<\dots<j_{i}\leq k} 
    \left(q\wedge(\smallwedge_{l=1}^{i}p_{j_{l}})\right)\left(x,u_{j_{1}},\dots,u_{j_{i}} \right).
\end{equation}
\end{theorem}
\begin{proof}
To prove Equation (\ref{eq:Adetmain}), we assume that the vectors $(u_{i})_{i=1,\dots,k}$ are
  linearly independent and we complete it with elements $(u_{i})_{i=k+1,\dots,n}$ to get 
  a basis in $V$, and assume that $\omega\in\Lambda^{n}(V)$ is a nonzero $n$-form.
We compute the determinant from the following equation
\begin{equation*}
\omega(Au_{1},\dots,Au_{n})=\omega(u_{1},\dots,u_{n})\det A.
\end{equation*}
A simple expansion of the expression
\begin{equation*}
\omega\left(u_{1}+\sum_{j=1}^{k}u_{j}p_{j}(u_{1}),\dots,
  u_{n}+\sum_{j=1}^{k}u_{j}p_{j}(u_{n})\right)
\end{equation*}
  gives Equation (\ref{eq:Adetmain}).

Let us denote the right hand side of Equation (\ref{eq:Ainvmain}) by $\varphi$, 
  define $m=\min(n-1,k)$, and assume that $m>1$, since the $m=1$ case is trivial.
If $x=Ay$ then we have the following equation.
\begin{align*}
\varphi=&q(y)+\sum_{a=1}^{k}q(u_{a})p_{a}(y)
  +\sum_{j=1}^{k}(q\wedge p_{j})(Ay,u_{j})\\
&+\sum_{i=2}^{m-1} \sum_{1\leq j_{1}<\dots<j_{i}\leq k}
  \left[q\wedge\left(\smallwedge_{l=1}^{i}p_{j_{l}}\right)\right](Ay,u_{j_{1}},\dots,u_{j_{i}})\\
& +\left[q\wedge\left(\smallwedge_{l=1}^{k}p_{l}\right)\right](Ay,u_{1},\dots,u_{k})
\end{align*}
We can expand the third term
\begin{equation*}
\sum_{j=1}^{k}(q\wedge p_{j})(Ay,u_{j})=q(y)\sum_{j=1}^{k}p_{j}(u_{j})
 -\sum_{j=1}^{k}q(u_{j})p_{j}(y)+\sum_{j=1}^{k}\sum_{\btop{c=1}{c\neq j}}^{k}
  p_{c}(y)(q\wedge p_{j})(u_{c},u_{j})
\end{equation*}
  the summand in the fourth term
\begin{align*}
&
\left[q\wedge\left(\smallwedge_{l=1}^{i}p_{j_{l}}\right)\right](Ay,u_{j_{1}},\dots,u_{j_{i}})=
  q(y)\left[q\wedge\left(\smallwedge_{l=1}^{i}p_{j_{l}}\right)\right]
  (u_{j_{1}},\dots,u_{j_{i}})+\\
&\qquad+\sum_{c=1}^{i}(-1)^{c}p_{j_{c}}(y)
  \left[q \wedge \left(\smallwedge_{\btop{l=1}{l\neq c}}^{i} p_{j_{l}}\right)\right]
(u_{j_{1}},\dots,u_{j_{i}})\\
&\qquad+\sum_{\btop{c=1}{c\notin\left\{j_{1},\dots,j_{i}\right\}}}^{k}p_{c}(y)
  \left[q\wedge\left(\smallwedge_{l=1}^{i}p_{j_{l}}\right)\right]  
  (u_{c},u_{j_{1}},\dots,u_{j_{i}})
\end{align*}
  and the fifth term
\begin{align*}
\left[q\wedge\left(\smallwedge_{l=1}^{k}p_{l}\right)\right](Ay,u_{1},\dots,u_{k})=&
  q(y)\left(\smallwedge_{l=1}^{k}p_{l}\right)(u_{1},\dots,u_{k})\\
&+ \sum_{c=1}^{k}(-1)^{c}p_{j_{c}}(y)
  \left[q \wedge \left(\smallwedge_{\btop{l=1}{l\neq c}}^{k} p_{j_{l}}\right)\right]
  (u_{1},\dots,u_{k}).
\end{align*}
Combining these terms we get
\begin{align*}
\varphi=&q(y)\left[1+\sum_{j=1}^{k}p_{j}(u_{j})
  +\sum_{i=2}^{m-1}\left(\smallwedge_{l=1}^{i}p_{j_{l}}\right)(u_{j_{1}},\dots,u_{j_{i}})
  +\left(\smallwedge_{l=1}^{k}p_{l}\right)(u_{1},\dots,u_{k}) \right]\\
&+\sum_{i=1}^{m-1}\sum_{1\leq j_{1}<\dots<j_{i}\leq k}
  \sum_{\btop{c=1}{c\notin\left\{j_{1},\dots,j_{i}\right\}}}^{k}p_{c}(y)
  \left[q\wedge\left(\smallwedge_{l=1}^{i}p_{j_{l}}\right)\right]
  (u_{c},u_{j_{1}},\dots,u_{j_{i}})\\
&+\sum_{i=2}^{m}\sum_{1\leq j_{1}<\dots<j_{i}\leq k}\sum_{c=1}^{i}
(-1)^{c}p_{j_{c}}(y)
  \left[q \wedge \left(\smallwedge_{\btop{l=1}{l\neq c}}^{i} p_{j_{l}}\right)\right]
  (u_{j_{1}},\dots,u_{j_{i}}).
\end{align*}
This can be rewritten as
\begin{align}
\varphi=q(y)\det A&\nonumber\\
\label{eq:Ainvproof1}
+\sum_{i=1}^{m-1}\Biggl\{&\sum_{1\leq j_{1}<\dots<j_{i}\leq k}
  \sum_{\btop{c=1}{c\notin\left\{j_{1},\dots,j_{i}\right\}}}^{k}p_{c}(y)
  \left[q\wedge\left(\smallwedge_{l=1}^{i}p_{j_{l}}\right)\right]
  (u_{c},u_{j_{1}},\dots,u_{j_{i}}) \\
+&\sum_{1\leq j_{1}<\dots<j_{i+1}\leq k}\sum_{c=1}^{i+1}
  (-1)^{c}p_{j_{c}}(y)
  \left[q \wedge \left(\smallwedge_{\btop{l=1}{l\neq c}}^{i+1} p_{j_{l}}\right)\right]
  (u_{j_{1}},\dots,u_{j_{i+1}})\Biggr\}. \nonumber
\end{align}
For a given $1\leq i\leq m-1$, we assume that $1\leq j_{1}<\dots<j_{i}\leq k$ are fixed.
Now we check how the term
\begin{equation}
\label{eq:Ainvproof2}
q\wedge\left(\smallwedge_{l=1}^{i}p_{j_{l}}\right)
\end{equation}
  occurs in the previous formula.
From the first summation in Equation (\ref{eq:Ainvproof1}), we have
\begin{equation}
\label{eq:Ainvproof3}
\sum_{\btop{c=1}{c\notin\left\{j_{1},\dots,j_{i}\right\}}}^{k}p_{c}(y)
  \left[q\wedge\left(\smallwedge_{l=1}^{i}p_{j_{l}}\right)\right]
  (u_{c},u_{j_{1}},\dots,u_{j_{i}}).
\end{equation}
If in the second summation we have the indices $1\leq j_{1}'<\dots<j_{i+1}'\leq k$,
  we get the term (\ref{eq:Ainvproof2}) if 
\begin{equation*}
(j_{1}',\dots,j_{b-1}',j_{b+1}',\dots,j_{i+1}')=(j_{1},\dots,j_{i})
\end{equation*}
  holds for a $b$ index.
This $b$ is the extra $j'$ index in the second summation.
To get the term (\ref{eq:Ainvproof2}) we have $c=b$, so the second summation is
\begin{equation*}
\sum_{\btop{b=1}{b\notin\left\{j_{1},\dots,j_{i}\right\}}}^{k}
  (-1)^{b}p_{j_{b}'}(y)
  \left[q \wedge \left(\smallwedge_{\btop{l=1}{l\neq b}}^{i+1} p_{j_{l}'}\right)\right]
  (u_{j_{1}'},\dots,u_{b},\dots u_{j_{i+1}'}).
\end{equation*}
This can be written as
\begin{equation}
\label{eq:Ainvproof4}
-\sum_{\btop{b=1}{b\notin\left\{j_{1},\dots,j_{i}\right\}}}^{k}p_{j_{b}}(y)
  \left[q\wedge\left(\smallwedge_{l=1}^{i}p_{j_{l}}\right)\right](u_{b},u_{j_{1}},\dots,u_{j_{i}}).
\end{equation}
Adding Equations (\ref{eq:Ainvproof3},\ref{eq:Ainvproof4}), we get zero, which means that
  the summands for all $1\leq i\leq m-1$ in Equation (\ref{eq:Ainvproof1}) are zero.
This proves that $\varphi=q(y)\det A$, which is the left hand side of the Equation
  (\ref{eq:Ainvmain}).
\end{proof}

We note, that the previous theorem can be proved by induction on $k$, but the detailed proof
  has approximately the same length. 
Now we can state our main result as a simple consequence of the previous theorem.

\begin{corollary}
Assume that $B:V\to V$ is an invertible map and consider the perturbated operator
\begin{equation}
B'=B+\sum_{i=1}^{k}v_{i}\otimes p_{i},
\end{equation}
  where $v_{i}\in V$ and $p_{i}\in V^{*}$ for all $i=1,\dots,k$.
Let us define $u_{i}=B^{-1}v_{i}$ (for all $i=1,\dots,k$) and
\begin{equation}
A=\id_{V}+\sum_{i=1}^{k} u_{i}\otimes p_{i}.
\end{equation}
If $\det A\neq 0$, then we have for the inverse of the perturbated operator
\begin{equation}
(B')^{-1}=\frac{1}{\det A}B^{-1}+\frac{1}{\det A}\left(\sum_{i=1}^{\min(n-1,k)}
  \sum_{1\leq j_{1}<j_{2}<\dots<j_{i}\leq k} 
  \bigosquare_{l=1}^{i}(u_{j_{l}},p_{j_{l}})\right)B^{-1},
\end{equation}
  where $n=\dim V$.
\end{corollary}
\begin{proof}
Since $B'=BA$, we use the formula $(B')^{-1}=A^{-1}B^{-1}$, where $A^{-1}$ is given by 
  Equation (\ref{eq:Ainv}).
\end{proof}

\section{Connection with the Taylor expansion}

In applications, we assume that the perturbation
\begin{equation}
B'=B+\sum_{i=1}^{k}v_{i}\otimes p_{i}=B+Q
\end{equation}
  is small in some sense with respect to $B$.
(In our framework there is no norm, so the word \textit{small} has just intuitive meaning here.)
If we take into account only $m$ ($m\in\mathbb{N}$) (or a less number of) products of 
  $(u_{i},p_{i})$, we get the $m$-th order approximation of $(B')^{-1}$, that is
\begin{equation*}
(B')^{-1}_{m}=\frac{1}{\sideset{}{_{m}}\det A}B^{-1}
  +\frac{1}{\sideset{}{_{m}}\det A}\left(\sum_{i=1}^{m}
  \sum_{1\leq j_{1}<j_{2}<\dots<j_{i}\leq k} 
  \bigosquare_{l=1}^{i}(u_{j_{l}},p_{j_{l}})\right)B^{-1},
\end{equation*}
  where
\begin{equation*}
\sideset{}{_{m}}\det A=1+\sum_{i=1}^{m} \alpha_{i},
  \quad
  \alpha_{i}=\sum_{1\leq j_{1}<j_{2}<\dots<j_{i}\leq k} 
      \left(\smallwedge_{l=1}^{i}p_{j_{l}}\right)\left(u_{j_{1}},\dots,u_{j_{i}} \right)
  \quad\mbox{and}\quad 
  u_{i}=B^{-1}v_{i}.
\end{equation*}
We have for the zeroth, the first and second order approximation of $(B')^{-1}$
\begin{align*}
&(B')^{-1}_{0}=B^{-1}\\
&(B')^{-1}_{1}=B^{-1}-\dfrac{1}{1+\alpha_{1}}B^{-1}QB^{-1}\\
&(B')^{-1}_{2}=B^{-1}-\frac{1+\alpha_{1}}{1+\alpha_{1}+\alpha_{2}}B^{-1}QB^{-1}+
  \frac{1}{1+\alpha_{1}+\alpha_{2}}B^{-1}QB^{-1}QB^{-1}
\end{align*}
  and in general
\begin{equation}
\label{eq:approx}
(B')^{-1}_{m}=B^{-1}+\dfrac{1}{\displaystyle 1+\sum_{i=1}^{m} \alpha_{i}}\times
  \sum_{i=1}^{m}(-1)^{i}\left(1+\sum_{j=1}^{m-i}\alpha_{j} \right)(B^{-1}Q)^{i}B^{-1}.
\end{equation}
It is obvious from the construction that $(B')^{-1}_{m}=(B')^{-1}$ if $m\geq k$.

Since the $i$-th derivative of the inversion function $\iota(B)=B^{-1}$ is
\begin{equation*}
\dint^{i}\iota(B)(Q^{\left[i\right]})=(-1)^{i}i!(B^{-1}Q)^{i}B^{-1},
\end{equation*}
  the $m$-th order Taylor expansion of $\iota$ is
\begin{equation}
\label{eq:Taylor}
T_{m}(B+Q)^{-1}=B^{-1}+\sum_{i=1}^{m}(-1)^{i}(B^{-1}Q)^{i}B^{-1}.
\end{equation}
It means that if we define $\alpha_{i}=0$ ($i=1,\dots,m$) then our approximation
 (Equation (\ref{eq:approx})) gives back the Taylor expansion.
However, these $\alpha_{i}$ parameters guarantee that our $m$-th order approximation
  gives the exact result if $m\geq\rank Q$, while the Taylor expansion gives just an
  approximation for every $m$.

\begin{remark}
Assume that $g$ is a metric on V, that is a bilinear symmetric map $g:V\times V\to\mathbb{R}$
  which is non-degenarate: for every $0\neq v\in V$ there exists a vector $u\in V$ 
  such that $g(u,v)\neq 0$.
Then for every $v\in V$ we have $g(v,\cdot)\in V^{*}$, and we can define an isomorphism
\begin{equation*}
\tilde{g}:V\to V^{*}\qquad v\mapsto g(v,\cdot).
\end{equation*}
Assume that $A:V\to V^{*}$ is a linear map and define $\Tilde{A}= \tilde{g}^{-1}A$ which is a 
  $V\to V$ linear map.
If $A$ is invertible and perturbated to $A'=A+W$, then the above-mentioned Theorem and 
  Corollary gives approximations $(\Tilde{A'})^{-1}_{m}$ for $(\Tilde{A'})^{-1}$ and we have
  the approximations $(\Tilde{A'})^{-1}_{m} \tilde{g}^{-1}$ for $(A')^{-1}$.
\end{remark}

\begin{remark}
If $\dim V\in\left\{2,\dots,10\right\}$ and the matrix $B$ is a random matrix, 
  $k\in\left\{2,\dots,15\right\}$ and the vectors $(v_{i})_{i=1,\dots,k}$ and covectors 
  $(p_{i})_{i=1,\dots,k}$ are random vectors, and $V$ is endowed with the Euclidean metric, then
  numerical simulations show, that the convergence of the given approximation is faster
  than the convergence of the Taylor expansion.
We conjecture that this numerical observation is true in general settings too.
\end{remark}

\end{document}